\documentclass[12pt]{article}
\usepackage{amsmath,amssymb,epsfig,subfigure,fancybox}
\usepackage{graphicx}
\usepackage{enumerate}
\usepackage{cite}
\newcommand{\lan}{\big\langle}
\newcommand{\ran}{\big\rangle}

\newcommand{\lbar}{\overline}
\newcommand{\wdh}{\widehat}
\newcommand{\wdt}{\widetilde}

\newcommand{\Dl}{\Delta}
\def\op{{\mathcal L}}

\newcommand{\F}{{\mathcal F}}
\newcommand{\e}{\varepsilon}
\newcommand{\rr}{{\Bbb R}}
\newcommand{\M}{{\mathcal M}}
\newcommand{\cd}{(\cdot)}

\def\para#1{\vskip .2\baselineskip\noindent{\bf #1}}
\def\qed{\strut\hfill $\Box$}
\def\eqdef {\stackrel {\rm def}{=}}

\newtheorem{thm}{Theorem}[section]

\newtheorem{lem}[thm]{Lemma}

\newtheorem{rem}[thm]{Remark}
\newtheorem{defn}[thm]{Definition}
\newtheorem{exm}[thm]{Example}

\newcommand{\thmref}[1]{Theorem~{\rm \ref{#1}}}
\newcommand{\lemref}[1]{Lemma~{\rm \ref{#1}}}

\def\al{\alpha}

\makeatletter \@addtoreset{equation}{section}

\newcommand{\beq}[1]{\begin{equation} \label{#1}}
\newcommand{\eeq}{\end{equation}}
\newcommand{\bed}{\begin{displaymath}}
\newcommand{\eed}{\end{displaymath}}
\newcommand{\bea}{\bed\begin{array}{rl}}
\newcommand{\eea}{\end{array}\eed}
\newcommand{\ad}{&\!\!\!\disp}
\newcommand{\aad}{&\disp}
\newcommand{\barray}{\begin{array}{ll}}
\newcommand{\earray}{\end{array}}
\def\({\left(}
\def\){\right)}
\def\disp{\displaystyle}

\begin{document}
\title{Mean-Variance Type Controls Involving a Hidden Markov Chain: Models and Numerical Approximation\thanks{This
research was supported in part by the National Science Foundation
under DMS-1207667.}}
\author{Zhixin Yang,\thanks{Department of Mathematics,
Wayne State University, Detroit, Michigan 48202,
zhixin.yang@wayne.edu} \and George Yin,\thanks{Department of
Mathematics, Wayne State University, Detroit, Michigan 48202,
gyin@math.wayne.edu}
\and
Qing Zhang\thanks{Department of Mathematics, The
University of Georgia, Athens, GA 30602, qingz@math.uga.edu.}}
\date{}
\maketitle

\begin{abstract}
Motivated by applications arising in networked systems,
this work examines controlled regime-switching systems that stem from
 a mean-variance formulation.
A main point is that
the switching process is a hidden Markov chain. An additional
piece of information, namely, a noisy observation of switching process corrupted by
white noise is available.
We focus on
 minimizing the variance subject to
a fixed terminal expectation. Using  the Wonham filter, we convert
the partially observed system
 to
a completely observable one first.
Since closed-form solutions are virtually impossible be obtained, a Markov chain approximation method is used to devise
a computational scheme. Convergence of the algorithm is obtained.
A numerical example is provided to demonstrate the results.

\vskip 0.2 true in
\para{Key Words.} Mean-variance control, numerical method, Wonham filter.
\end{abstract}

\newpage

\setlength{\baselineskip}{0.25in}
\section{Introduction}
Using a switching diffusion model,
in our recent work \cite{YYWZ}, three potential applications in
platoon controls were outlined based on mean-variance controls. The first
concerns the longitudinal
 inter-vehicle distance control. To increase
 highway utility, it is desirable
 to reduce the total length of a platoon,
 resulting in reducing inter-vehicle distances.
 This strategy, however,
 increases the risk of collision in the presence of
 vehicle traffic uncertainties. To minimize the risk
  with desired inter-vehicle distance
 can be mathematically modeled as a mean-variance
 optimization problem.
 The second one is communication resource
 allocation of bandwidths for vehicle to vehicle (V2V)
 communications. For a given maximum throughput
 of a platoon communication system,
 the communication system operator
 must find a way to assign this resource
 to different V2V channels, which
  may also be formulated as a mean-variance control problem.
 The third one
 is the platoon fuel consumption that is
 total vehicle fuel consumptions within the platoon.
 Due to variations in vehicle sizes and speeds,
 each vehicle's fuel consumption is a controlled
 random process.
 Tradeoff between a platoon's team acceleration/maneuver
 capability and fuel consumption can be summarized in
 a desired platoon fuel consumption rate.
 Assigning fuels to different vehicles result
 in coordination of vehicle operations modeled
 by subsystem fuel rate dynamics. This problem
  may also be formulated  as a mean-variance control problem.

To capture the underlying dynamics of these problems,
it is natural to model the underlying system as diffusions coupled
by a finite-state Markov chain. For example, in the first case
of applications, the Markov chain may represent
external and macro states including
traffic states (road condition, overall congestions),
weather conditions (major thunder/snow storms), etc.
These macro states are observable with some noise.

This paper extends the mean-variance methods
to incorporate possible hidden Markov chains and to apply the results to
network control problems. In particular,
the underlying system is modeled as a controlled switching diffusion
modulated  by a finite-state Markov chain representing the system modes.
The state of the Markov chain is observable with additive white noise.
Given the target expectation of the state variable at the terminal time,
the objective is to minimize the variance at the terminal.
We use the mean-variance approach to treat the problem
and aim at developing feasible numerical methods
for solutions of the associated control problems.

Ever since the classical Nobel prize winning mean-variance portfolio
selection models for a single period was established by Markowitz in
\cite{M}, there has been much effort devoted to studying modern
portfolio theory in finance. Extensions toward different directions
have been pursued (for example, \cite{S,P}). Continuous-time mean-variance hedging problems were also examined; see \cite{D} among
others, in which hedging contingent claims in incomplete markets
problem was considered and optimal dynamic strategies were obtained
with the help of projection theorem. In the traditional set up, the
tradeoff between the risk and return is usually implicit, which
makes the investment decision much less intuitive. Zhou and Li
\cite{ZL} introduced an alternative methods to deal with the
mean-variance problems in continuous time, which embedded the
original problem into a tractable auxiliary problem, following  Li
and Ng's paper \cite{LN} for the multi-period model. They were able
to solve the auxiliary problem explicitly by linear quadratic theory
with the help of backward stochastic differential equations; see the
linear quadratic control problems with indefinite control weights in
\cite{LZ} and also \cite{JZ} and references therein. Recently, much
attention has been drawn to modeling controlled systems with
random environment and other factors that cannot be completely
captured by a simple diffusion model. In this connection, a set of
diffusions with regime switching appears to be suitable for the
problem.
Regime-switching models have been used in options pricing
\cite{DZ}, stock selling rules \cite{Z}, and mean-variance models
\cite{ZY} and \cite{Y2}.
The regime-switching models have also been considered in our work
\cite{YYWZ} using a two-time-scale formulation.

In connection with network control problems,
  while the current paper concentrates on the formulation and
  numerical methods. Detailed treatment of the specific platoon
  applications will be considered in a separate paper.
In our formulation, the coefficients of the systems are modulated by
a Markov chain. In contrast to many models in the literature, the
Markov chain is hidden, i.e., it is not completely observable.
In this paper, we consider the case that
a function of the chain with additive noise is observable.
In networked systems, such measurement can be
obtained with the addition of a sensor.

The underlying problem is a stochastic control problems with partial observation.
To resolve the problem,
we resort to Wonham filter to estimate the state. Then the original system
is converted into
a completely observable one.
In stochastic control literature,
a suboptimal filter for linear systems
with hidden Markov switching coefficients was considered
in \cite{DB} in connection with a quadratic cost control problem.
In this paper, we formulate the problem as a Markov modulated
mean-variance control problem with partial information.
Under our formulation, it is difficult to obtain a closed-form solution
 in contrast to \cite{ZY}. We need to resort to numerical
 algorithms.
We
use the Markov chain approximation methods of Kushner and
Dupuis \cite{KD} to develop
numerical algorithms.
Different from \cite{SY} and \cite{LZY},  the variance
is control dependent. In view of this, extra care must be taken to
address such control dependence.
The main purpose of this paper is to develop numerical methods
for the partially observed mean-variance control problem.
Applications in networked systems including implementation issues will
be considered elsewhere.

Starting from the partially observed control problems, our contributions
of this paper include:
\begin{itemize}
\item[{(1)}]  We use Wonham filtering techniques to convert
the problem into a completely observable system.

\item[{(2)}] We develop numerical approximation techniques based on
the Markov chain approximation schemes.
Although Markov chain approximation techniques have been used extensively in
various stochastic systems, the work on combination of such a methods with partial observed
control systems seems to be scarce to the best of our knowledge.
Different from the existing work in the literature, we use
Markov chain approximation for the diffusion component and use a direct
discretization for the Wonham filter.

\item[{(3)}] We use weak convergence methods to obtain the convergence of the algorithms.
A feature that is different from the existing work is that
 in the martingale problem formulation, the states include a component that comes from
Wonham filtering.
\end{itemize}

The rest of the paper is arranged as follows. Section \ref{sec:for}
presents the
problem formulation. Section
\ref{sec:mar} introduces the Markov
chain approximation methods. Section \ref{sec:opt} deals with the
approximation of the optimal controls. In Section \ref{sec:con}, we establish the
convergence of the algorithm. Section \ref{sec:num} gives one numerical
example for illustration; also included are some further remarks to conclude the paper.

\section{Formulation}\label{sec:for}
This section presents the formulation of the problem. We begin with
notation and assumptions. Given a probability space  $(\Omega,
\mathcal{F}, P)$ in which there are $w_1(t)$, a standard $d$
dimensional Brownian motion with $w_1(t)=(w^1_1(t),w^2_1(t),\ldots,
w^d_1(t))'$ where $z'$ denotes the transpose of $z$, and  a
continuous-time finite states Markov chain $\alpha(t)$ that is
independent of $w_1(t)$ and that takes values in $\M=\{1,2,\ldots,
m\}$ with generator $Q=(q_{ij})_{m\times m}$. We consider such a
networked system that there are $d+1$ nodes
in which one of the nodes follows the stochastic ODE%
\beq{flow1}\barray dx_1(t)\ad=r(t,\al(t))x_1(t)dt, \ \ t\in[s,T]\\
x_1(s)\ad=x_1,\earray\eeq where $r(t,i)\ge0$ for $i=1,2,\ldots,m$ is
the increase rate corresponding to different regimes in the network
systems. The flows of other $d$ nodes $x_l(t)$, $l=2,3,\ldots,d+1$
satisfy the system
of SDEs%
\beq{flowd}\barray
dx_l(t)\ad=x_l(t)b_l(t,\al(t))dt+x_l(t)\bar{\sigma}_l(t,\al(t))dw_1(t)\\
\ad=x_l(t)b_l(t,\al(t))dt+x_l(t)\bar{\sigma}_l(t,\al(t))dw_1(t),t\in
[s,T]\\
x_l(s)\ad =x_l,\earray\eeq where for each $i$, $b_l(t,i)$ is the
increase rate process and
$\bar{\sigma}_l(t,i)=(\bar{\sigma}_{l1}(t,i),\ldots,\bar{\sigma}_{ld}(t,i))$
is the volatility for the $l$th node. In our framework, instead of
having full information of the Markov chain, we can only observe it
in white noise. That is, we observe $y(t)$, whose dynamics is given
by \beq{obs} \barray
 dy(t)=g(\alpha(t))dt+\sigma_0 dw_2(t),\\
 y(s)=0,
 \earray
  \eeq
 where $\sigma_0>0$ and $w_2\cd$ is
a standard scalar Brownian motion, where $w_2(\cdot)$, $w_1(\cdot)$,
and $\alpha(\cdot)$ are independent. Moreover, the
initial data
$p(s)=p=(p^1, p^2,\ldots,p^m)$ in which $p^i=p^i(s)=P(\al(s)=i)$ is
given for $1\le i\le m$. By distributing $N_l(t)$ shares of flows to
$l$th node at time $t$ and denoting the total flows for the whole
networked system as $x(t)$ we have
 $$x(t)=\sum^{d+1}_{l=1}N_l(t)x_l(t), t\ge s.$$   Therefore, the
dynamics of
 $x(t)$ are given as
\beq{sys}
 \barray
 dx(t)\ad=\sum^{d+1}_{l=1}N_l(t)dx_l(t)\\
 \ad=[r(t,\al(t))N_1(t)x_1(t)+\sum^{d+1}_{l=2}b_l(t,\al(t))N_l(t)x_l(t)]dt\\
 \aad \hspace*{1in}
 +\sum^{d+1}_{l=2}N_l(t)x_l(t)\sum^d_{j=1}\bar{\sigma}_{lj}(t,\al(t))dw^j_1(t)\\
 \ad=[r(t,\al(t))x(t)+\sum^{d+1}_{l=2}(b_{l}(t,\al(t))-r(t,\al(t)))u_l(t)]dt
 +\sum^{d+1}_{l=2}\sum^d_{j=1}\bar{\sigma}_{lj}(t,\al(t))u_l(t)dw^j_1(t)\\
 \ad=[x(t)r(t,\alpha(t))+B(t,\alpha(t))u(t)]dt+u'(t)\bar{\sigma}(
 t,\alpha(t))dw_1(t),\\
 x(s)\ad=\sum^{d+1}_{l=1}N_l(s)x_l(s)=x, 
\earray \eeq in which $u(t)=(u_2(t),\ldots,u_{d+1}(t))'$ and
$u_l(t)=N_l (t)x_l(t)$ for $l=2,\ldots, d+1$ is the actual flow of
the network system for the $l$th node and
$u_1(t)=x(t)-\sum^{d+1}_{l=2}u_l(t)$ is the actual flow of the
networked system for the first node, and
\bea \ad B(t,\al(t))=(b_2(t,\al(t))-r(t,\al(t)),\ldots,b_{d+1}(t,\al(t))-r(t,\al(t))),\\
\ad
\bar{\sigma}(t,\al(t))=(\bar{\sigma}_1(t,\al(t)),\bar{\sigma}_2(t,\al(t)),\ldots,\bar{\sigma}_d(t,\al(t)))'=(\bar{\sigma}_{lj}(t,\al(t)))_{d\times
d}
.\eea We define $\F_t=\sigma\{w_1(\wdt s),y(\wdt s),x(s): s \le \wdt
s\le t\}$. Our objective is to find an $\F_t$ admissible control
$u(\cdot)$ in a compact set $\mathcal{U}$ under the constraint that
the expected terminal flow value is $Ex(T)=\kappa$ for some given
$\kappa\in \mathbb{R}$, so that the risk measured by the variance of
terminal flow is minimized. Specifically, we have the following goal
\beq{obj} \barray \ad \text{ min }
J(s,x,p,u(\cdot)):=E[x(T)-\kappa]^2\\
\ad \text{ subject to } Ex(T)=\kappa. \earray \eeq To handle the
constraint part in problem \eqref{obj}, we apply the Lagrange
multiplier technique and thus get unconstrained problem (see,
e.g.,\cite{ZL}) with multiplier $\lambda$:
\beq{obj2}%
\barray%
 \ad \text{ min }
J(s,x,p,u(\cdot),\lambda):=E[x(T)+\lambda-\kappa]^2-\lambda^2\\
\ad \text{ subject to } (x(\cdot),u(\cdot)) \text{ admissible. }
\earray \eeq
A pair $(\sqrt{{\rm Var}\;( x(T))},\kappa)\in \rr^2$ corresponding to the optimal control,
if exists, is called
an {\it efficient point}.
The set of all the efficient points is called the {\it efficient frontier}.

 Note that one of the striking feature of our model is that we
have no access to the value of Markov chain at a given time $t$,
which makes the problem more difficult than\cite{ZY}. Let
$p(t)=(p^1(t),\ldots,p^m(t))\in \rr^{1\times m}$ with
$p^i(t)=P(\alpha(t)=i|{\cal F}^y(t))$ for $i=1,2,\ldots,m,$ with
$\mathcal{F}^y(t)=\sigma\{y(\wdt s): s \le \wdt s\le t\}$. It was
shown in Wonham \cite{Wo} that this conditional probability satisfies
the following system of stochastic differential equations
\beq{won}\barray
dp^i(t)\ad=\sum_{j=1}^mq^{ji}p^j(t)dt+\frac{1}{\sigma_0}p^i(t)
(g(i)-\lbar\alpha(t))d\widehat w_2(t),\\
 p^i(s)\ad=p^i,\earray\eeq
where $\lbar\alpha(t)=\sum_{i=1}^mg(i)p^i(t)$ and $\widehat w_2(t)$
is the innovation process. It is easy to see that $\widehat
w_2(\cdot)$ is independent of $w_1(\cdot)$.

\begin{rem}
{\rm
Note that in connection with portfolio optimization, the additional
observation process $y(t)$ can be related to non-public (insider) information.
Insider information is often corrupted by noise and may reveal the
direction of the underlying security prices.
}
\end{rem}

\begin{rem}
{\rm In \cite{LZY}, a much simpler model was considered in
connection with an asset allocation problem. In particular, the
diffusion gain $\sigma$ is independent of $\al(t)$. This makes it
possible to convert the original system into a completely observable
one with the help of Wonham filter. Nevertheless, under our
framework, the dependence on $\al(t)$ in $\sigma$ is crucial and the
corresponding nonlinear filter is of infinity dimensional. In view
of this, we can only turn to approximation schemes.
}
\end{rem}

With the help of Wonham filter, given the independence conditions,
we can find the best estimator for $r(t,\alpha(t))$,
$B(t,\alpha(t))$, and $\bar{\sigma}(t,\alpha(t))$ in the sense of
least mean square prediction error and transform the partial
observable system into completely observable system given as below:
\bea
 dx(t) = [{\wdh{r(t,\alpha (t))}}x(t) + \wdh{B(t,\alpha (t))}u(t)]dt + u'(t)\wdh{\bar{\sigma} (t,\alpha
 (t))}dw_1(t),\eea where
 \beq{eq3}%
\barray%
\wdh {r(t,\al(t))}\ad\eqdef \sum^m_{i=1}r(t,i)p^i(t)\in \rr^1,\\
\wdh {B(t,\al(t))}\ad\eqdef
(\sum^m_{i=1}(b_2(t,i)-r(t,i))p^i(t),\ldots,\sum^m_{i=1}(b_{d+1}(t,i)-r(t,i))p^i(t))\in\rr^{1\times d},\\
\wdh {\bar{\sigma}(t,\al(t))}\ad\eqdef (\sum^m_{i=1}
\bar{\sigma}_{lj}(t,i)p^i(t))_{d\times d}
.\earray%
\eeq%
Note that $u'(t)\wdh {\bar{\sigma}(t,\al(t))}$ is an $\rr^{1\times
d}$ row vector which is defined as \bea \disp u'(t)\wdh
{\bar{\sigma}(t,\al(t))}\ad =\sigma(x(t),p(t),u(t))\\
\ad =(\sigma_1(x(t),p(t),u(t)),\sigma_2(x(t),p(t),u(t)),\ldots,\sigma_d(x(t),p(t),u(t))).\eea
In this way, by putting the two components $p(t)$ and $x(t)$
together, we get $$(x(t),p(t))=(x(t),p^1(t),...,p^m(t)),$$ a
completely observable system  whose dynamics are as follows
\beq{eq2}%
 \barray
 dx(t)\ad= [\sum^m_{i=1}r(t,i)p^i(t)x(t) + \sum^{d+1}_{l=2}\sum^m_{i=1}(b_l(t,i)-r(t,i))p^i(t)u_l(t)]dt\\
 \aad\ \  + \sum^{d+1}_{l=2}\sum^d_{j=1}\sum^m_{i=1}u_l(t)\bar{\sigma}_{lj}(t,i)p^i(t)dw^j_1(t) \\
\ad=b(x(t),p(t),u(t))dt+\sigma(x(t),p(t),u(t))dw_1(t)\\
dp^i(t)\ad=\sum_{j=1}^mq^{ji}p^j(t)dt+\frac
{1}{\sigma_0}p^i(t)(g(i)-\lbar\alpha(t))d\widehat
w_2(t), \text{
for }i=\{1,\ldots,m\}\\
x(s)\ad=x,\quad  p^i(s)=p^i.\\
 \earray
\eeq%

To proceed, for an arbitrary $r\in \mathcal{U}$ and
$\phi(\cdot,\cdot,\cdot)\in C^{1,2,2}(\rr)$, we first define the
differential operator $\op^r$ by%
\beq{op}\barray \op^r\phi(s,x,p)\ad=\frac{\partial \phi}{\partial
s}+ \frac{\partial \phi}{\partial
x}b(x,p,r)+\frac{1}{2}\frac{\partial^2\phi}{\partial
x^2}[\sigma(x,p,r)\sigma'(x,p,r)]\\
\aad\ \ +\sum^m_{i=1}\frac{\partial \phi}{\partial
p^i}\sum_{j=1}^mq^{ji}p^j
+\frac{1}{2}\sum^m_{i=1}\frac{\partial^2\phi}{\partial
(p^i)^2}\frac{1}{\sigma^2_0}[p^i(g(i)-\lbar\alpha)]^2.\earray \eeq%

Let $W(s,x,p,u)$ be the objective function and let $E^u_{s,x,p}$
denote the expectation of functionals on $[s,T]$ conditioned on
$x(s)=x, p(s)=p$ and the admissible control $u=u\cd$.
\beq{cost}\barray
W(s,x,p,u)=E^u_{s,x,p}(x(T)+\lambda-k)^2-\lambda^2\earray\eeq
and $V(s,x,p)$ be the value function \beq{vf}\barray
V(s,x,p)=\inf_{u\in \mathcal{U}}W(s,x,p,u).\earray\eeq

The value function is a solution of the following system of HJB
equation%
\beq{eq5}\barray
   \inf_{r\in\mathcal{U}} \op^r V(s,x,p)=0,\earray\eeq
with boundary condition $V(T,x,p)=(x+\lambda-\kappa)^2-\lambda^2$.

We have successfully  converted
 an optimal control problem with partial observations to a problem with full
 observation. Nevertheless, the problem has not been completely solved.
Due to the high nonlinearity and complexity, a closed-form solution
of the optimal control problem is virtually impossible to obtain. As
a viable alternative, we use the Markov chain approximation
techniques \cite{KD} to construct numerical schemes to approximate
the optimal strategies and the optimal values. Different from the
standard numerical scheme, we construct a discrete-time controlled Markov chain
to approximate the diffusions of the $x\cd$ process. For the
 Wonham filtering equation, we approximate the solution by discretizing it directly.
In fact, to implement the Wonham filter, we take logarithmic
transformation to discretize the resulting equation.

\section{Discrete-time Approximation Scheme}\label{sec:mar}

In this section, we deal with the numerical algorithms
for the two components system. First, for the second component
$p^i(t)$, numerical experiments and simulations show that
discretizing the stochastic differential equation about $p^i(t)$
directly could produce undesirable results (such as producing
a non-probability vector and numerically unstable etc.)
due to white noise perturbations.
It may produce a non-probability result. To overcome this
difficulty, we use the idea in  \cite[Section 8.4]{YinZ} and
transform the dynamic system of $p^i(t)$, then design a numerical
procedure for the
transformed system. Let $v^i(t)=\log p^i(t)$ and apply
 the It\^{o}'s rule lead to the following dynamics to obtain
\beq{eq22}%
\barray%
dv^i(t) \ad=[\sum^m_{j=1}
q^{ji}\frac{p^j(t)}{p^i(t)}-\frac{1}{2\sigma^2_0}(g(i)-\bar{\al}(t))^2]dt
+\frac{1}{\sigma_0}[g(i)-\bar{\al}(t)]d\widehat
w_2(t),\\
v^i(s)\ad=\log(p^i).
\earray%
\eeq%

By choosing the constant step size $h_2>0$ for time variable we can
discrete \eqref{eq22} as follows:
\beq{dis}\barray%
v^{h_2,i}_{n+1}\ad=v^{h_2,i}_n+h_2[\sum^m_{j=1
}q^{ji}\frac{p^{h_2,j}_n}{p^{h_2,i}_n}-\frac{1}{2\sigma^2_0}(g(i)-
\bar{\al}^{h_2}_n)^2]+\sqrt{h_2}\frac{1}{\sigma_0}(g(i)-\bar{\al}^{h_2}_n)\e_n,\\
v^{h_2,i}_0\ad=\log(p^{i}),\\
p^{h_2,i}_{n+1}\ad=\exp(v^{h_2,i}_{n+1}),\\
p^{h_2,i}_0\ad=p^i,\earray\eeq where
$\bar{\al}^{h_2}_n=\sum^m_{i=1}g(i)p^{h_2,i}_n$ and $\{\e_n\}$ is a
sequence of i.i.d. random variables satisfying $E\e_n=0$, $E\e^2_n=1$,
and $E|\e_n|^{2+\gamma}<\infty$ for some $\gamma>0$ with
 $$\e_n=\frac{
\wdh  {w}_2((n+1)h_2)-\wdh  {w}_2(nh_2)}{\sqrt{h_2}}.$$  Note that
$p^{h_2,i}_n$ appeared as the denominator in \eqref{dis} and we have
focused on the case that $p^{h_2,i}_n$ stays away from $0$. A
modification can be made to take into consideration the case of
$p^{h_2,i}_n=0$. In that case, we can choose a fixed yet arbitrarily
large positive real number $M$ and use the idea of penalization to
construct the approximation as below:
\beq{dis1}\barray%
v^{h_2,i}_{n+1}\ad=v^{h_2,i}_n+h_2\{[\sum^m_{j=1
}q^{ji}\frac{p^{h_2,j}_n}{p^{h_2,i}_n}-\frac{1}{2\sigma^2_0}(g(i)-\bar{\al}^{h_2}_n)^2]I_{\{p^{h_2,i}_n\ge
 e^{-M}\}}-MI_{\{p^{h_2,i}_n<e^{-M}\}}\}\\
\aad\quad +\sqrt{h_2}\frac{1}{\sigma_0}(g(i)-\bar{\al}^{h_2}_n)\e_n,\\
v^{h_2,i}_0\ad=\log(p^i),\\
p^{h_2,i}_{n+1}\ad=\exp(v^{h_2,i}_{n+1}),\\
p^{h_2,i}_0\ad=p^i.\earray\eeq

In what follows, we construct a discrete-time finite state Markov
chain to approximate the controlled diffusion process, $x(t)$. Given
that in our model, we have both time variable $t$ and state variable
$p(t)$ and $x(t)$ involved. Our construction of Markov chain needs
to take care of time and state variables as follows. Let $h_1>0$ be
a discretizatioin parameter for state variables, and recall that
$h_2>0$ is the step size for time variable. Let $N_{h_2}=(T-s)/h_2$
be an integer and define
$S_{h_1}=\{x:x=kh_1,k=0,\pm1,\pm2,\ldots\}$. We use $u^{h_1,h_2}_n$
to denote the random variable that is the control action for the
chain at discrete time $n$. Let $u^{h_1,h_2}=(u^{h_1,h_2}_0,
u^{h_1,h_2}_1,\ldots)$ denote the sequence of $\mathcal{U}$-valued
random variables which are the control actions at time $0,1,\ldots$
and  $p^{h_2}=(p^{h_2}_0, p^{h_2}_1,\ldots)$ are the corresponding
posterior probability in which
$p^{h_2}_n=(p^{h_2,1}_n,p^{h_2,2}_n,\ldots,p^{h_2,m}_n)$. We define
the difference $\Delta
\xi^{h_1,h_2}_n=\xi^{h_1,h_2}_{n+1}-\xi^{h_1,h_2}_n$ and let
$E^{h_1,h_2,r}_{x,p,n}$, $Var^{h_1,h_2,r}_{x,p,n}$ denote the
conditional expectation and variance given
$\{\xi^{h_1,h_2}_k,u^{h_1,h_2}_k,p^{h_2}_k, k\le n,
\xi^{h_1,h_2}_n=x,p^{h_2}_n=p, u^{h_1,h_2}_n=r\}$. By stating that
$\{\xi^{h_1,h_2}_n,n<\infty\}$ is a controlled discrete-time Markov
chain on a  discrete time  state space $S_{h_1}$ with transition
probabilities from state $x$ to another state $y$, denoted by
$p^{h_1,h_2}((x,y)|r,p)$, we mean that the transition probabilities
are functions of a control variable $r$ and posterior probability
$p$. The sequence $\{{\xi^{h_1,h_2}_n,n<\infty}\}$ is said to be
locally consistent with \eqref{eq2}, if it satisfies
\beq{eq3.3}%
\barray%
E^{h_1,h_2,r}_{x,p,n}\Dl\xi^{h_1,h_2}_n=b(x,p,r)h_2+o(h_2),\\
V^{h_1,h_2,r}_{x,p,n}\Dl\xi^{h_1,h_2}_n=\sigma(x,p,r)\sigma'(x,p,r)h_2+o(h_2),\\
\sup_n|\Dl\xi^{h_1,h_2}_n|\to 0,\text{  as  }h_1,h_2\to 0.
\earray%
\eeq%
Let $\mathcal{U}^{h_1,h_2}$ denote the collection of ordinary
controls, which is determined by a sequence of such measurable
functions $F^{h_1,h_2}_n(\cdot)$  that
$u^{h_1,h_2}_n=F^{h_1,h_2}_n(\xi^{h_1,h_2}_k, p^{h_2}_k, k\leq n,
u^{h_1,h_2}_k, k<n)$. We say that $u^{h_1,h_2}$ is admissible for
the chain if $u^{h_1,h_2}_n$ are $\mathcal{U}$ valued random
variables and the Markov property continues to hold under the use of
the sequence $\{{u^{h_1,h_2}_n}\}$, namely, \bea \ad
P\{\xi^{h_1,h_2}_{n+1}=y|\xi^{h_1,h_2}_k,u^{h_1,h_2}_k,
p^{h_2}_k,k\le
n\}\\
\ad=P\{\xi^{h_1,h_2}_{n+1}=y|\xi^{h_1,h_2}_n,u^{h_1,h_2}_n,p^{h_2}_n\}=p^{h_1,h_2}((\xi^{h_1,h_2}_n,y)|u^{h_1,h_2}_n,p^{h_2}_n).\eea
With the approximating Markov chain given above, we can approximate
the objective function defined in \eqref{cost} by%
\beq{eq3.1}%
W^{h_1,h_2}(s,x,p,u^{h_1,h_2})=E^{u^{h_1,h_2}}_{s,x,p}(\xi^{h_1,h_2}_{N_{h_2}}+\lambda-k)^2-\lambda^2.%
\eeq%
Here, $E^{u^{h_1,h_2}}_{s,x,p}$ denote the expectation given that
$\xi^{h_1,h_2}_0=x, p^{h_2}_0=p$ and that  an admissible control
sequence $u^{h_1,h_2}=\{{u^{h_1,h_2}_n, n<\infty}\}$ is used.
 Now we need the approximating Markov chain constructed above
 satisfying local consistency, which is one of the
necessary conditions for weak convergence. To find a reasonable
Markov chain that is locally consistent, we first suppose that
control space has a unique admissible control
$u^{h_1,h_2}\in\mathcal{U}^{h_1,h_2}$, so that we can drop inf in
\eqref{eq5}. We discrete \eqref{op} by the following finite
difference method using step-size $h_1>0$ for state variable and
$h_2>0$ for time variable as mentioned above. 
 \beq{d0}\barray V(t,x,p)\to
V^{h_1,h_2}(t,x,p);\earray\eeq For the derivative with respect to
the time variable, we use \beq{dt}\barray V_t(t,x,p)\to
\frac{V^{h_1,h_2}(t+h_2,x,p)-V^{h_1,h_2}(t,x,p)}{h_2};\earray\eeq
For the first derivative with respect to $x$, we use one-side difference method%
\beq{di}%
\barray%
V_{x}(t,x,p)\ad \to\left\{ {\begin{array}{*{20}c}
   \frac{V^{h_1,h_2}(t+h_2,x+h_1,p)-V^{h_1,h_2}(t+h_2,x,p)}{h_1} & \text{ for } {b(x,p,r)\ge 0} \\
   \frac{V^{h_1,h_2}(t+h_2,x,p)-V^{h_1,h_2}(t+h_2,x-h_1,p)}{h_1} & \text{ for }{b(x,p,r)<0} .
\end{array}} \right.\\
\earray \eeq
For the second derivative with respect to $x$, we have standard difference method%
 \beq{dii}\barray
V_{xx}(t,x,p)\ad\to\frac{V^{h_1,h_2}(t+h_2,x+h_1,p)+V^{h_1,h_2}(t+h_2,x-h_1,p)-2V^{h_1,h_2}(t+h_2,x,p)}{h^2_1}.
\earray%
\eeq%
For the first and second derivative with respect to posterior
probability, we also have the similar expression as above. Let
$V^{h_1,h_2}(t,x,p)$ denote the solution to the finite difference
equation with $x$ and $p$ be an integral multiplier of $h_1$ and $nh_2<T$.
Plugging all the necessary expressions into \eqref{eq5} and combining
the like terms and multiplying all terms by $h_2$ yield the
following expression:
\beq{dif} \barray \ad\hspace*{-0.5in} V^{h_1,h_2}(nh_2,x,p) \\
\ad\hspace*{-0.3in} =V^{h_1,h_2}(nh_2+h_2,x,p)[1-\frac{|b(x,p,r)|h_2}{h_1}-\frac{h_2\sigma(x,p,r)\sigma'(x,p,r)}{h^2_1}]\\
\aad\hspace*{-0.3in}+V^{h_1,h_2}(nh_2+h_2,x+h_1,p)\frac{\sigma(x,p,r)\sigma'(x,p,r)h_2+2h_1h_2b^+(x,p,r)}{2h^2_1}\\
\aad\hspace*{-0.3in}+V^{h_1,h_2}(nh_2+h_2,x-h_1,p)\frac{\sigma(x,p,r)\sigma'(x,p,r)h_2+2h_1h_2b^-(x,p,r)}{2h^2_1}\\
\aad\hspace*{-0.3in}+\sum^m_{i=1}V^{h_1,h_2}(nh_2+h_2,x,p^i+h_1)\frac{\frac{1}{\sigma^2_0}[p^i(g(i)-\lbar\alpha)]^2h_2+2h_1(\sum_{j=1}^mq^{ji}p^j)^+h_2}{2h^2_1}\\
\aad\hspace*{-0.3in}+\sum^m_{i=1}V^{h_1,h_2}(nh_2+h_2,x,p^i-h_1)\frac{\frac{1}{\sigma^2_0}[p^i(g(i)-\lbar\alpha)]^2h_2+2h_1(\sum_{j=1}^mq^{ji}p^j)^-h_2}{2h^2_1}\\
\aad\hspace*{-0.3in}+\sum^m_{i=1}V^{h_1,h_2}(nh_2+h_2,x,p^i)[-\frac{\frac{1}{\sigma^2_0}
[p^i(g(i)-\lbar\alpha)]^2h_2}{h^2_1}-\frac{h_2|\sum^m_{j=1}q^{ji}p^j|}{h_1}],
\earray\eeq where $b^+(x,p,r)$, $(\sum_{j=1}^mq^{ji}p^j)^+$ and
$b^-(x,p,r)$, $(\sum_{j=1}^mq^{ji}p^j)^-$ are positive and negative
parts of  $b(x,p,r)$ and $\sum_{j=1}^mq^{ji}p^j$, respectively. Note
the sum of the coefficient of the first three line in the above
equation is unity. By choosing proper $h_1$ and $h_2$, we can
reasonably assume that the coefficient
$$1-\frac{|b(x,p,r)|h_2}{h_1}-\frac{h_2\sigma(x,p,r)\sigma'(x,p,r)}{h^2_1}$$
of term  $V^{h_1,h_2}(nh_2+h_2,x,p)$ is in $[0,1]$. Therefore, the
coefficients can be regarded as the transition function of a Markov
chain. We define the transition probability in the following way,
\beq{tran}%
\barray%
\disp
p^{h_1,h_2}((nh_2,nh_2+h_2))|x,p,r)=1-\frac{|b(x,p,r)|h_2}{h_1}-\frac{h_2\sigma(x,p,r)\sigma'(x,p,r)}{h^2_1}\\
\disp
p^{h_1,h_2}((nh_2,x),(nh_2+h_2,x+h_1)|p,r)=\frac{\sigma(x,p,r)\sigma'(x,p,r)h_2+2h_1h_2b^+(x,p,r)}{2h^2_1}\\
\disp
p^{h_1,h_2}((nh_2,x),(nh_2+h_2,x-h_1)|p,r)=\frac{\sigma(x,p,r)\sigma'(x,p,r)h_2+2h_1h_2b^-(x,p,r)}{2h^2_1}\\
\earray%
\eeq%
Theoretically, we can find approximation of $V(s,x,p)$ in \eqref{vf}
by using \eqref{eq3.1}
and%
\beq{val} V^{h_1,h_2}(s,x,p)=\inf_{u^{h_1,h_2}\in
\mathcal{U}^{h_1,h_2}}W^{h_1,h_2}(s,x,p,u^{h_1,h_2}).\eeq
Practically, with the transition probability defined as above, we
can compute $V^{h_1,h_2}(s,x,p)$ by the following iteration method
\beq{tran1}\barray \ad\hspace*{-0.5in} V^{h_1,h_2}(nh_2,x,p)\\
\ad \hspace*{-0.3in}=
p^{h_1,h_2}((nh_2,x)(nh_2+h_2,x+h_1)|p,r)V^{h_1,h_2}(nh_2+h_2,x+h_1,p)\\
\aad\hspace*{-0.3in}+p^{h_1,h_2}((nh_2,x),(nh_2+h_2,x-h_1)|p,r)V^{h_1,h_2}(nh_2+h_2,x-h_1,p)\\
\aad\hspace*{-0.3in}+p^{h_1,h_2}((nh_2, nh_2+h_2)|x,p,r)V^{h_1,h_2}(nh_2+h_2,x,p)\\
\aad\hspace*{-0.3in}+\sum^m_{i=1}V^{h_1,h_2}(nh_2+h_2,x,p^i+h_1)\frac{\frac{1}{\sigma^2_0}[p^i(g(i)-\lbar\alpha)]^2h_2+2h_1(\sum_{j=1}^mq^{ji}p^j)^+h_2}{2h^2_1}\\
\aad\hspace*{-0.3in}+\sum^m_{i=1}V^{h_1,h_2}(nh_2+h_2,x,p^i-h_1)\frac{\frac{1}{\sigma^2_0}[p^i(g(i)-\lbar\alpha)]^2h_2+2h_1(\sum_{j=1}^mq^{ji}p^j)^-h_2}{2h^2_1}\\
\aad\hspace*{-0.3in}+\sum^m_{i=1}V^{h_1,h_2}(nh_2+h_2,x,p^i)[-\frac{\frac{1}{\sigma^2_0}[p^i(g(i)-\lbar\alpha)]^2h_2}{h^2_1}-\frac{h_2|\sum^m_{j=1}q^{ji}p^j|}{h_1}].
\earray\eeq
Note that we used local transitions here, we can avoid
the problem of ``numerical noise" or ``numerical viscosity" in this
way,  which appears  in non-local transitions case, and is even more
serious in higher dimension, see\cite{K1} for more details.  We can
show that the Markov chain $\{\xi^{h_1,h_2}_n, n<\infty\}$ with
transition probability $p^{h_1,h_2}(\cdot)$ defined in \eqref{tran}
is locally consistent with \eqref{eq2} by verifying the following
equations:
\beq{eq3.4}%
\barray%
\ad E^{h_1,h_2,r}_{x,p,n}\Dl
\xi^{h_1,h_2}_n\\
\aad\ \ =h_1\left(\frac{\sigma(x,p,r)\sigma'(x,p,r)h_2+2h_1h_2b^+(x,p,r)}{2h^2_1}\right)\\
\aad\ \ \ \ -h_1\left(\frac{\sigma(x,p,r)\sigma'(x,p,r)h_2+2h_1h_2b^-(x,p,r)}{2h^2_1}\right)\\
\aad\ \ =b(x,p,r)h_2,\\[2ex]
\ad V^{h_1,h_2,r}_{x,p,n}\Dl\xi^{h_1,h_2}_n\\
\aad\ \
=h^2_1\left(\frac{\sigma(x,p,r)\sigma'(x,p,r)h_2+2h_1h_2b^+(x,p,r)}{2h^2_1}\right)
\\
\aad\ \ \ \
+h^2_1\left(\frac{\sigma(x,p,r)\sigma'(x,p,r)h_2+2h_1h_2b^-(x,p,r)}{2h^2_1}\right)\\
\aad\ \ =\sigma(x,p,r)\sigma'(x,p,r)h_2+O(h_1h_2).\\[2ex]%
\earray%
\eeq%
\section{Approximation of Optimal Controls}\label{sec:opt}
\subsection{Relaxed Control and Martingale Measure}
Note the fact that the sequence of ordinary control constructed in
Markov chain approximation scheme may not converge in a traditional
sense due to the issue of closure. That is, a bounded sequence
$\xi^{h_1,h_2}_n$ with ordinary controls $u^{h_1,h_2}_n$ would not
necessarily have a subsequence which converges to a limit process
which is a solution to the equation driven by a desirable  ordinary
control. The use of the relaxed control gives us an alternative to
obtain and characterize the weak limit appropriately. Although the
usage of relaxed control enlarges the control space of the problem,
it does not alter the infimum of the objective function. We first
give the definition of
relaxed control as follows.%
\begin{defn}
For the $\sigma$-algebra $\mathcal{B}(\mathcal{U})\ \  and\ \
\mathcal{B}(\mathcal{U}\times [s,T])$ of Borel subsets of
$\mathcal{U}$ and $\mathcal{U}\times[s,T]$, an admissible relaxed
control or simply a relaxed control $m(\cdot)$ is a measure on
$\mathcal {B}(\mathcal {U}\times[s,T])$ such that
$m(\mathcal{U}\times[s,t])=t-s$ for all $t\in[s,T]$.
\end{defn}
For notional simplicity, for any $B\in \mathcal{B}(\mathcal{U})$, we
write $m(B\times[s,T])$ as $m(B,T-s)$. Since
$m(\mathcal{U},t-s)=t-s$ for all $t\in[s,T]$ and $m(B,\cdot)$ is
nondecreasing, it is absolutely continuous. Hence the derivative
$\dot{m}(B,t)=m_t(B)$ exists almost everywhere for each $B$. We can
further define the relaxed control representation $m(\cdot)$ of
$u(\cdot)$ by\beq{rc}\barray
 m_t(B)=I_{\{u(t)\in B\}}\earray
 \text{ for  any } B\in \mathcal{B}(\mathcal{U}).\eeq
 Therefore, we can represent any ordinary admissible control
 $u(\cdot)$ as a relaxed control by using the
 definition $m_t(dr)=I_{u(t)}(r)dr$, where $I_u(r)$ is the
 indicator function concentrated at the point $u=r$.  Thus, the
 measure-valued derivative $m_t(\cdot)$ of the relaxed control
 representation of $u(t)$ is a measure which is concentrated at the
 point $u(t)$.
 For each $t$, $m_t(\cdot)$ is a measure on
$ \mathcal{B}(\mathcal{U})$ satisfying $m_t(\mathcal{U})=1$
 and $m(A)=\int_{\mathcal{U}\times[s,T]}I_{\{(r,t)\in A\}}m_t(dr)dt$
for all $A\in \mathcal{B}(\mathcal{U}\times[s,T])$, i.e.,
$m(drdt)=m_t(dr)dt$. %

On the other hand, note that we have control in the diffusion gain.
The similar problem arises even with the introduction of relaxed
control. Therefore, we need to borrow the idea of martingale measure
to allow the desired closure and at the same time keep the same
infimum for the objective function. We say that $M(\cdot)$ is a
measure-value $\F_t$ martingale with values $M(B,t)$ if $M(B,\cdot)$
is an $\F_t$ martingale for each $B\in\cal{U}$, and for each $t$, the
following hold: $\sup_{B\in\cal{U}}EM^2(B,t)<\infty$, $M(A\cup B,
t)=M(A,t)+M(B,t)$ w.p.1. for all disjoint $A,B\in \cal{U}$, and
$EM^2(B_n,t)\to 0$ if $B_n\to \emptyset$. $M(\cdot)$ is said to be
continuous if each $M(B,\cdot)$ is. We say that $M(\cdot)$ is orthogonal
if $M(A,\cdot)$, $M(B,\cdot)$ is an $\mathcal{F}_t$ martingale
whenever $A\cap B=\emptyset$. If $M(\cdot)$, $\bar{M}(\cdot)$ are
$\mathcal{F}_t$ martingale measures and $M(A,\cdot)$,
$\bar{M}(B,\cdot)$ are $\mathcal{F}_t$ martingales for all Borel set
$A,B$, then $M(\cdot)$ and $\bar{M}(\cdot)$ are said to be strongly
orthogonal. Let $M(\cdot)=(M_1(\cdot),\ldots,M_d(\cdot))'$, a vector
valued martingale measure, we impose the following
conditions.
\begin{itemize}
\item [{(A1)}] $M(\cdot)=(M_1(\cdot), \ldots,M_d(\cdot))'$ is square
integrable and continuous, each component is orthogonal, and the
pairs are strongly orthogonal.

Under this assumption, there are measure-valued random processes
$m_i(\cdot)$ such that the quadratic variation processes satisfies,
for each $t$ and $B\in\mathcal{U}$ $$\lan M_i(A,\cdot), M_j(B,
\cdot)\ran (t)=\delta_{ij}m_i(A\cap B,t).$$

\item[{(A2)}] The $m_i$'s  do not depend on $i$,
so $m_i\cd=m\cd$, and $m(U,t)=t$ for all $t$.
\end{itemize}
With the use of relaxed control representation, the operator of the
controlled diffusion is given by \beq{op1}\barray
\mathcal{L}^mf(s,x,p)\ad=f_s+\int f_xb(x,p,c)m_t(dc)+\frac{1}{2}\int
f_{xx}\sigma(x,p,r)\sigma'(x,p,r)m_t(dc)\\
\aad+\sum^m_{i=1}
f_{p^i}\sum_{j=1}^mq^{ji}p^j+\frac{1}{2}\sum^m_{i=1}f_{p^ip^i}
\frac{1}{\sigma^2_0}[p^i(g(i)-\lbar\alpha)]^2\\
\ad=\int \mathcal{L}^rf(s,x,p)m_t(dc).\earray\eeq

Let there be a continuous process $(x(\cdot),p(\cdot))$ and a
measure $m(\cdot)$ satisfying assumption (A1) and (A2)  such that
for each bounded and smooth function $f(\cdot,\cdot,\cdot)$,
$$f(t,x(t),p(t))-f(s,x,p)-\int\int
\mathcal{L}^rf(z,x(z),p(z))m_z(dc)dz=Q_f(t)$$ is an $ \mathcal{\wdt
F}_t$ martingale, where $\mathcal{\wdt F}_t$ measures
$\{x(z),p(z),m_z(\cdot),s\le z\le t\}$. Then
$(x(\cdot),p(\cdot), m(\cdot))$ solves the martingale problem with
operator $\mathcal{L}^r$ and there is a martingale measure $M(\cdot)$
with quadratic variation $m(\cdot)I$ satisfying assumption (A1) and
(A2) such that
\beq{rx}\barray%
 x(t)\ad=x+\int^t_s\int_\mathcal{U}b(x(z),p(z),c)m_z(dc)dz+\int^t_s
\int_\mathcal{U}\sigma(x(z),p(z),c)M(dc,dz)\\
p^i(t)\ad=\int^t_s\sum_{j=1}^mq^{ji}p^j(z)dz+\int^t_s
\frac{1}{\sigma_0}[p^i(z)(g(i)-\lbar\alpha(z))]d\wdh {w}_2(z),\text{
for }i=\{1,\ldots,m\},\earray\eeq
where  $$\sigma(x(z),p(z),c)=(\sigma_1(x(z),p(z),c),\ldots,
\sigma_d(x(z),p(z),c))\in \rr^{1\times d}.$$
 Equation \eqref{rx} represents our control system. In the next
 section, we work on approximation of $(x(t),p(t),M(t),m(t))$. We
 say that $(M(\cdot),m(\cdot))$ is an admissible relaxed control
  for \eqref{rx} if $(A1)$ and $(A2)$ hold and $\lan M(\cdot)\ran =m(\cdot)I$. To
 proceed, we first suppose that
 \begin{itemize}
 \item[{(A3)}] $b(\cdot,\cdot,\cdot)$, $\sigma(\cdot,\cdot,\cdot)$ are
 continuous, $b(\cdot, p, c)$, $\sigma(\cdot,p,c)$ are
 Lipschitz continuous uniformly in $p$, $c$ and bounded.
 \item[{(A4)}] $\sigma(x,p,r)=(\sigma_1(x,p,r),\ldots,\sigma_d(x,p,r))>0$
 \end{itemize}

\subsection{Approximation of $(x(t),p(t),M(t), m(t))$}

Using $E^{h_1,h_2}_n$ to denote the conditional expectation given
$\{\xi^{h_1,h_2}_k, p^{h_2}_k, u^{h_1,h_2}_k,k\le n\}$. Define
$R^{{h_1,h_2}}_n=(\Dl\xi^{h_1,h_2}_n-E^{h_1,h_2}_n\Dl\xi^{h_1,h_2}_n)$.
By local consistency, we have
$$\xi^{h_1,h_2}_{n+1}=\xi^{h_1,h_2}_n+b(\xi^{h_1,h_2}_n, p^{h_2}_n,u^{h_1,h_2}_n)h_2+R^{{h_1,h_2}}_n,$$
where $cov^{h_1,h_2}_nR^{{h_1,h_2}}_n=a(\xi^{h_1,h_2}_n,
p^{h_2}_n,u^{h_1,h_2}_n)=\sum^d_{j=1}\sigma^2_j(\xi^{h_1,h_2}_n,
p^{h_2}_n,u^{h_1,h_2}_n)h_2+O(h_1h_2)$. Note that we can decompose
$a(\xi^{h_1,h_2}_n,
p^{h_2}_n,u^{h_1,h_2}_n)=P^{h_1,h_2}_n(D^{h_1,h_2}_n)^2(P^{h_1,h_2}_n)'$,
in which $P^{h_1,h_2}_n=(\frac{1}{\sqrt{d}},\cdots,
\frac{1}{\sqrt{d}})\in \rr^{1\times d}$ and $D^{h_1,h_2}_n$ is
diagonal
$$D^{h_1,h_2}_n=\{\sqrt{d}\sigma_1(\xi^{h_1,h_2}_n,
p^{h_2}_n,u^{h_1,h_2}_n),\sqrt{d}\sigma_2(\xi^{h_1,h_2}_n,
p^{h_2}_n,u^{h_1,h_2}_n),\cdots, \sqrt{d}\sigma_d(\xi^{h_1,h_2}_n,
p^{h_2}_n,u^{h_1,h_2}_n) \}\in \rr^{d \times d},$$ then we can
represent $R^{{h_1,h_2}}_n$ in terms of Brownian motion defined as
$$\Delta
w^{h_1,h_2}_n=(D^{h_1,h_2}_n)^{-1}(P^{h_1,h_2}_n)'R^{h_1,h_2}_n.$$
In this way, $R^{h_1,h_2}_n=\sigma(\xi^{h_1,h_2}_n,
p^{h_2}_n,u^{h_1,h_2}_n)\Dl w^{h_1,h_2}_n+\e^{h_1,h_2}_n$
(see \cite[Section10.4.1]{KD} for details). We can thus
represent $\xi^{h_1,h_2}_{n+1}$ as \beq{dx}\barray
\xi^{h_1,h_2}_{n+1}=\xi^{h_1,h_2}_n+b(\xi^{h_1,h_2}_n,
p^{h_2}_n,u^{h_1,h_2}_n)h_2 +\sigma(\xi^{h_1,h_2}_n,
p^{h_2}_n,u^{h_1,h_2}_n)\Delta
w^{h_1,h_2}_n+\e^{h_1,h_2}_n.\earray\eeq To take care of the control
part, let $\{C^{h_1,h_2}_l,l\le k_{h_1,h_2}\}$ be a finite
 partition of $\mathcal{U}$ such that the diameters of $C^{h_1,h_2}_l\to 0$
 as $h_1,h_2\to 0$. Let $c_l \in C^{h_1,h_2}_l$. Define the random variable
$$\Delta w^{h_1,h_2}_{l,n}=\Delta
w^{h_1,h_2}_nI_{\{u^{h_1,h_2}_n=c_l\}}+\Delta
\psi^{h_1,h_2}_{l,n}I_{\{u^{h_1,h_2}_n \neq c_l\}}.$$  Then we have%
 \beq{fd}\barray
 \ad\xi^{h_1,h_2}_{n+1}=\xi^{h_1,h_2}_n+b(\xi^{h_1,h_2}_n, p^{h_2}_n, u^{h_1,h_2}_n)h_2
 +\sum_{l=1}\sigma(\xi^{h_1,h_2}_n, p^{h_2}_n, u^{h_1,h_2}_n)I_{\{u^{h_1,h_2}_n=c_l\}}\Delta w^{h_1,h_2}_{l,n}+\e^{h_1,h_2}_n,\\
\ad m^{h_1,h_2}_n(c_l)=I_{\{u^{h_1,h_2}_n=c_l\}}. \earray\eeq In
order to approximate the continuous time process $(x(t),p(t),M(t),
m(t))$, we use continuous-time interpolation. We define the
piecewise constant interpolations by
 \beq{ite}\barray
\ad \xi^{h_1,h_2}(t)=\xi^{h_1,h_2}_n, \quad p^{h_2}(t)=p^{{h_2}}_n,
\bar{\al}^{h_1,h_2}(t)=\sum^m_{i=1}g(i)p^{h_2}_n,
u^{h_1,h_2}(t)=u^{h_1,h_2}_n,\\
\ad z^{h_2}(t)=n, w^{h_1,h_2}_l(t)=\sum^{z^{h_2}(t)-1}_{k=0}\Delta
w^{h_1,h_2}_{l,k}, \e ^{h_1,h_2}(t)=\e^{h_1,h_2}_n, \quad
\text{for}\quad t\in[nh_2,(n+1)h_2).\earray\eeq
 Define relaxed representation $m^{h_1,h_2}(\cdot)$ of $u^{h_1,h_2}(\cdot)$ by
$m^{h_1,h_2}_t(B)=I_{\{u^{h_1,h_2}(t)\in B\}}$
 for any $B\in \mathcal{B}(\mathcal{U})$. $m^{h_1,h_2}(dc,dt)=m^{h_1,h_2}_t(dc)dt$ and
  $m^{h_1,h_2}_t(\cdot)=m^{h_1,h_2}_n(\cdot)$ for $t\in[nh_2,nh_2+h_2)$.
Here a
 sequence $m^{h_1,h_2}_n(\cdot)$ of measure-valued random variables is an
 admissible relaxed control if $m^{h_1,h_2}_n(\mathcal{U})=1$ and $${P\{\xi^{h_1,h_2}_{n+1}=
 y|\xi^{h_1,h_2}_i, p^{h_2}_i, m^{h_1,h_2}_i,i\le n\}}=\int p^{h_1,h_2}(\xi^{h_1,h_2}_n,y|p^{h_2}_n,c)m^{h_1,h_2}_n(dc).$$
 For $c_l\in C^{h_1,h_2}_l$, $\{M(C^{h_1,h_2}_l,\cdot),l\le
k_{h_1,h_2}\}$ are orthogonal continuous
   martingale with $\lan M(C^{h_1,h_2}_l,\cdot)\ran =m(C^{h_1,h_2}_l,\cdot)$. There are mutually
 independent $d$ dimensional standard Wiener process
 $w^{h_1,h_2}_l(\cdot), l\le
 k_{h_1,h_2}$ such that
 \beq{mr}\barray
 M(C^{h_1,h_2}_l,t)=\int^t_s(m_z(C^{h_1,h_2}_l))^{\frac{1}{2}}dw^{h_1,h_2}_l(z).\earray\eeq
Let $M^{h_1,h_2}(\cdot)$ and $m^{h_1,h_2}(\cdot)$ be the
restrictions of the measures of $M(\cdot)$ and $m(\cdot)$,
respectively, to the sets$\{C^{h_1,h_2}_l,l\le k_{h_1,h_2}\}$. The
following lemma demonstrate the fact that we can approximate
$(x(t),p(t),M(t), m(t))$ by a quadruple satisfying
\beq{appr}\barray
 \xi^{h_1,h_2}(t)\ad=x+\int^t_s\int_\mathcal{U}b(\xi^{h_1,h_2}(z),p^{h_2}(z),c)m^{h_1,h_2}_z(dc)dz\\
 \aad+\int^t_s\int_\mathcal{U}\sigma(\xi^{h_1,h_2}(z),p^{h_2}(z),c)M^{h_1,h_2}(dc,dz)+\e^{h_1,h_2}(t)\\
\ad=x+\int^t_s\sum_lb(\xi^{h_1,h_2}(z),p^{h_2}(z),c_l)m_z(C^{h_1,h_2}_l)dz\\
\aad+\int^t_s\sum_l\sigma(\xi^{h_1,h_2}(z),p^{h_2}(z),c_l)(m_z(C^{h_1,h_2}_l))^\frac{1}{2}dw^{h_1,h_2}_l(z)
+\e^{h_1,h_2}(t),
 \earray\eeq where $m^{h_1,h_2}(\cdot)$ is a piecewise constant and takes
 finitely many values and $M^{h_1,h_2}(\cdot)$ is represented in
 terms of a finite number of Wiener process. The idea is similar to the method used in \cite[Theorem 8.1]{K}, we
omit the detail here for brevity.
 \begin{lem}\label{ap}
 Assume $(A1)-(A4)$ and satisfying
 \eqref{appr}, then $$(\xi^{h_1,h_2}(\cdot), p^{h_2}(\cdot), m^{h_1,h_2}(\cdot), M^{h_1,h_2}(\cdot))\Rightarrow (x(\cdot), p(\cdot), m(\cdot), M(\cdot)).$$
Also, $W(s,x,p,m^{h_1,h_2})\to W(s,x,p,m)$ and we can suppose that
$m^{h_1,h_2}(\cdot)$ is piecewise constant further.
 \end{lem}

 Let $\mathcal{F}^{h_1,h_2}_t$ denote the
$\sigma$-algebra that measures at least
 \beq{si}\barray
 \{\xi^{h_1,h_2}(z),p^{h_2}(z),m^{h_1,h_2}_z(\cdot), M^{h_1,h_2}(\cdot), w^{h_1,h_2}_l(z),1\le l\le k_{h_1,h_2} ,s \le z\le t\}.\earray
 \eeq
Using $\Gamma^{h_1,h_2}$ to denote the set of admissible relaxed
control $m^{h_1,h_2}(\cdot)$ with respect to\\
$\{w^{h_1,h_2}_l(\cdot),p^{h_2}(\cdot),l\le k_{h_1,h_2}\}$ such that
$m^{h_1,h_2}_t(\cdot)$
 is a fixed probability measure in the interval $[nh_2,(n+1)h_2)$.  With the notation of relaxed control
 given above, we can write 
 \eqref{eq3.1} and value
 function \eqref{val} as
\beq{3eq3.1}%
W^{h_1,h_2}(s,x,p,m^{h_1,h_2})=E^{m^{h_1,h_2}}_{s,x,p}(\xi^{h_1,h_2}(T)+\lambda-k)^2-\lambda^2.%
\eeq%
\beq{2val}V^{h_1,h_2}(s,x,p)=\inf_{m^{h_1,h_2}\in\Gamma^{h_1,h_2}}W^{h_1,h_2}(s,x,p,m^{h_1,h_2}).\eeq

Note also that \eqref{cost} can be written in terms of the relaxed control:
\beq{cost2}\barray
W(s,x,p,m)=E^m_{s,x,p}(x(T)+\lambda-k)^2-\lambda^2.\earray\eeq

\section{Convergence}\label{sec:con}
Let
$(\xi^{h_1,h_2}(\cdot),p^{h_2}(\cdot),m^{h_1,h_2}(\cdot),M^{h_1,h_2}(\cdot))$
be a solution of \eqref{appr}, where $M^{h_1,h_2}(\cdot)$ is a martingale
measure with respect to the filtration $\mathcal{F}^{h_1,h_2}_t$,
with quadratic variation process $m^{h_1,h_2}(\cdot)$. Then we can
proceed to obtain the convergence of the algorithm next.

\begin{thm}\label{1}
Under Assumption {\rm(A1)-(A5)}. Let the approximating chain
$\{\xi^{h_1,h_2}_n, n<\infty \}$ be constructed with transition
probability defined in \eqref{tran}, and $p^{h_2}_n$ is approximated
by \eqref{dis}. Let $\{u^{h_1,h_2}_n, n<\infty\}$ be a sequence of
admissible controls, $\xi^{h_1,h_2}(\cdot)$ and $p^{h_2}(\cdot)$ be
the continuous time interpolation defined in \eqref{ite},
$m^{h_1,h_2}(\cdot)$ be the relaxed control representation of
$u^{h_1,h_2}(\cdot)$ $($continuous time interpolation of
$u^{h_1,h_2}_n)$. Then
$\{\xi^{h_1,h_2}(\cdot),p^{h_2}(\cdot),m^{h_1,h_2}(\cdot)\}$ is
tight. Denoting the limit of a weakly convergent subsequence by
$\{x(\cdot),p(\cdot), m(\cdot)\}$, there exists a martingale measure
$M(\cdot)$, with respect to $\{\mathcal{F}_t,t\ge s\}$, and with
quadratic variation process $m(\cdot)$ such that \eqref{rx} is
satisfied.
\end{thm}

\para{Proof.}
Note that $m^{h_1,h_2}(\cdot)$ is
tight due to the compactness of the relaxed control under the weak topology.
Since $(\xi^{h_1,h_2}(\cdot),p^{h_2}(\cdot)) \in \rr^{m+1}$,
the tightness of $p^{h_2}(\cdot)$ can be obtained as in \cite[Theorem
8.15]{YinZ}. Therefore, we just need to
take care that of $\xi^{h_1,h_2}(\cdot)$ in the following part.  For
the tightness of $\xi^{h_1,h_2}(\cdot)$, by assumption {(A1)}, for
$s\le t\le T$, \beq{sm}\barray
E^{m^{h_1,h_2}}_{s,x,p}|\xi^{h_1,h_2}(t)-x|^2\ad
=E^{m^{h_1,h_2}}_{s,x,p}|\int^t_s\int_\mathcal{U}
b(\xi^{h_1,h_2}(z),p^{h_2}(z),c)m^{h_1,h_2}_z(dc)dz\\
\aad\ \ \
+\int^t_s\int_\mathcal{U}\sigma(\xi^{h_1,h_2}(z),p^{h_2}(z),c)M^{h_1,h_2}(dc,dz)
\\
\aad \ \ \hfill +\e^{h_1,h_2}(t)|^2\\
\ad \le Kt^2+Kt+\e^{h_1,h_2}(t).\earray\eeq Here $K$ is a generic
positive constant whose value may be different in different context.
Similarly, we can guarantee
$E^{m^{h_1,h_2}}_{s,x,p}|\xi^{h_1,h_2}(t+\delta)-\xi^{h_1,h_2}(t)|^2=O(\delta)+\e^{h_1,h_2}(\delta)$
as $\delta\to 0$. Therefore, the tightness of
$\xi^{h_1,h_2}(\cdot)$ follows. By the compactness of set
$\mathcal{U}$, we can see that $M^{h_1,h_2}(\cdot)$ is also tight.
In view of the tightness, we can extract a
weakly convergent subsequence, and denote its limit by
$\{x(\cdot),p(\cdot), m(\cdot),M(\cdot)\}$. We next  show that the
limit is the solution of SDE driven by $(p(\cdot),m(\cdot),
M(\cdot))$.

For $\delta>0$ and any process $\nu(\cdot)$ define the process
$\nu^\delta(\cdot)$ by $\nu^\delta(t)=\nu(n\delta )$ for
$t\in[n\delta,n\delta+\delta)$. Then by the tightness of
$\xi^{h_1,h_2}(\cdot)$ and $p^{h_2}(\cdot)$, \eqref{appr} can be
rewritten as
 \beq{yd}\barray\disp
\xi^{h_1,h_2}(t)\ad=x+\int^t_s\int_\mathcal{U}
b(\xi^{h_1,h_2}(z),p^{h_2}(z),c)m^{h_1,h_2}_z(dc)dz\\
\aad\ \ +\int^t_s\int_\mathcal{U}
\sigma(\xi^{{h_1,h_2},\delta}(z),p^{h_2,\delta}(z),c)M^{h_1,h_2}(dc,dz)+\e^{{h_1,h_2},\delta}(t),\earray
\eeq %
where $\lim_{\delta\to 0}\lim\sup_{{h_1,h_2}\to
0}E|\e^{{h_1,h_2},\delta}(t)|\to 0$.\\ We further assume that the
probability space is chosen as required by Skorohod representation.
Therefore, we can assume the sequence
$\{\xi^{h_1,h_2}(\cdot),p^{h_2}(\cdot),m^{h_1,h_2}(\cdot),
M^{h_1,h_2}(\cdot)\}$ converges to $(x(\cdot),p(\cdot), m(\cdot),
M(\cdot))$ w.p.1 with a little bit abuse of notation.

Taking limit as $h_1\to 0$ and $h_2\to 0$, the convergence of
$\{\xi^{h_1,h_2}(\cdot),p^{h_2}(\cdot),m^{h_1,h_2}(\cdot),
M^{h_1,h_2}(\cdot)\}$ to its limit w.p.1 implies that
$$E|\int^t_s\int_\mathcal{U}b(\xi^{h_1,h_2}(z),p^{h_2}(z),c)m^{h_1,h_2}_z(dc)dz
-\int^t_s\int_\mathcal{U}b(x(z),p(z),c)m^{h_1,h_2}_z(dc)dz|\to 0,$$
uniformly in $t$. Also, recall that $m^{h_1,h_2}(\cdot)\to m(\cdot)$
in the ``compact weak" topology if and only if
$$\int^t_s\int_\mathcal{U}\phi(c,z)m^{h_1,h_2}(dc,dz)\to \int^t_s\int_\mathcal{U}
\phi(c,z)m(dc,dz).$$ for any continuous and bounded function
$\phi(\cdot)$ with compact support. Thus, weak convergence and
Skorohod representation imply that \beq{bp}\barray\disp
\int^t_s\int_\mathcal{U}b(x(z),p(z),c)m^{h_1,h_2}_z(dc)dz\to
\int^t_s\int_\mathcal{U}b(x(z),p(z),c)m_z(dc)dz \text{ as }
h_1,h_2\to 0,\earray\eeq uniformly in $t$ on any bounded interval
w.p.1.

Recall that $M^{h_1,h_2}(\cdot)$ is a martingale measure with
quadratic variation process $m^{h_1,h_2}(\cdot)$. Due to the fact
that $\xi^{h_1,h_2,\delta}(\cdot)$ and $p^{h_2,\delta}(\cdot)$ are
piecewise constant functions, following from the probability one
convergence, we have

\beq{dp}\barray\disp
\int^t_s\int_\mathcal{U}\sigma(\xi^{{h_1,h_2},\delta}(z),p^{h_2,\delta}(z),c)M^{h_1,h_2}(dc,
dz)\to \int^t_s\int_\mathcal{U}
\sigma(x^\delta(z),p^{\delta}(z),c)M^{h_1,h_2}(dc,dz).\earray\eeq
Recall that recall that $M^{h_1,h_2}(\cdot)\to M(\cdot)$ in the
``compact weak" topology if and only if
$\int^t_s\int_\mathcal{U}f(c,z)M^{h_1,h_2}(dc,dz)\to
\int^t_s\int_\mathcal{U}f(c,z)M(dc,dz)\text{ as } h_1,h_2\to 0$ for
each bounded and continuous function $f(\cdot)$, we have
$$\int^t_s \int_\mathcal{U}\sigma(x^\delta(z),p^{\delta}(z),c)M^{h_1,h_2}(dc,dz)\to
\int^t_s
\int_\mathcal{U}\sigma(x^\delta(z),p^{\delta}(z),c)M(dc,dz),$$
uniformly
 in $t$ on any bounded interval w.p.1; see \cite[pp. 352]{KD}. %
Combining the above results, we have \beq{fyd}\barray \disp
x(t)=x+\int^t_s\int_\mathcal{U}b(x(z),p(z),c)m(dc,dz)+\int^t_s\int_\mathcal{U}
\sigma(x^\delta(z),p^{\delta}(z),c)M(dc,dz)+\e^\delta(t).\earray\eeq
Where $\lim_{\delta\to 0}E|\e^\delta(t)|=0$. Taking limit of the
above equation as $\delta \to 0$ yields \eqref{rx}.\qed

\begin{thm}
  Under assumptions {\rm(A1)-(A5)}, $V^{h_1,h_2}(s,x, p)$ and $V(s,x,
  p)$ are value functions defined in \eqref{2val} and
  \eqref{vf} respectively, we have\beq{lim} \barray V^{h_1,h_2}(s,x,p)\to\
  V(s,x,p), \text{ as } h_1\to 0, h_2\to 0.\earray\eeq
\end{thm}

\para{Proof.} For each $h_1,h_2$, let $\wdh {m}^{h_1,h_2}$ be an
optimal relaxed control for $\{x^{h_1,h_2}(\cdot),p^{h_2}(\cdot)\}$.
i.e.
$$V^{h_1,h_2}(s,x, p)=W^{h_1,h_2}(s,x, p,\wdh {m}^{h_1,h_2})=\inf_{m^{h_1,h_2}
\in\Gamma^{h_1,h_2}}W^{h_1,h_2}(s,x, p,m^{h_1,h_2})$$
Choose a subsequence $\{\wdt h_1,\wdt h_2\}$ of $\{h_1,h_2\}$ such that
$$\lim \inf_{{h_1,h_2}\to 0} V^{h_1,h_2}(s,x, p)=\lim_{{\wdt h_1,\wdt h_2}\to 0}
V^{{\wdt h_1,\wdt h_2}}(s,x, p)
=\lim_{{\wdt h_1,\wdt h_2}\to 0}W^{\wdt h_1,\wdt h_2}(s,x, p,\wdh {m}^{\wdt h_1,\wdt h_2}).$$
Note that
we can assume that $\{\xi^{\wdt h_1,\wdt h_2}(\cdot),p^{\wdt h_2}(\cdot),\wdh
{m}^{\wdt h_1,\wdt h_2}(\cdot),\wdh {M}^{\wdt h_1,\wdt h_2}(\cdot)\}$ converges
weakly to $\{x(\cdot),p(\cdot),m(\cdot),M(\cdot)\}$. Otherwise, take
a subsequence of $\{\wdt h_1,\wdt h_2\}$ to assume its weak limit.
\thmref{1}, Skorohod representation and dominance convergence
theorem imply that as  $\wdt h_1,\wdt h_2\to 0$
$$E^{\wdh {m}^{\wdt h_1,\wdt h_2}}_{s,x,p}(\xi^{\wdt h_1,\wdt h_2}(T)+\lambda-k)^2-\lambda^2\to
E^m_{s,x,p}(x(T)+\lambda-k)^2-\lambda^2.$$ So
$$W^{\wdt h_1,\wdt h_2}(s,x,p, \wdh {m}^{\wdt h_1,\wdt h_2})\to W(s,x,p,m)\ge
V(s,x,p).$$ It follows that $$\lim\inf_{h_1,h_2\to 0}
V^{h_1,h_2}(s,x,p)\ge V(s,x,p)$$ Next, we need to show
$\lim\sup_{{h_1,h_2}\to 0}V^{h_1,h_2}(s,x,p)\le V(s,x,p)$ to
complete the proof.
 Given any $\rho>0$,
there is a $\delta>0$,
with the help
of \lemref{ap}, we are able to approximate any such quadruple
$(x(t), p(t), m(t), M(t))$  by a quadruple satisfying
 \bea
 x^\delta(t)\ad=x+\int^t_s\int_\mathcal{U}b(x^\delta(z),p^\delta(z),c)m^\delta_z(dc)dz+\int^t_s
\int_\mathcal{U}\sigma(x^\delta(z),p^\delta(z),c)M^\delta(dc,dz),\eea
where $m^\delta(\cdot)$ is  piecewise constant and takes finitely
many values and $M^\delta(\cdot)$ is represented in terms of a
finite number of $d$-dimensional Wiener process such that
for
the  optimization problem
with \eqref{rx}
and \eqref{cost2} under the constraints that the control are
concentrated on the points $c_1,c_2,\ldots, c_N$ for all $t$. They
take on one value $c_j$ on each interval
$[\iota\delta,\iota\delta+\delta), \iota=0,1,\ldots.$ Let $\wdh
{u}^\rho(\cdot)$ be the optimal control and $\wdh {m}^\rho(\cdot)$
be its relaxed control representation, and let $(\wdh
{x}^\rho(\cdot), \wdh{p}^{\rho}(\cdot))$ be the associated solution
process. Since $\wdh {m}^\rho(\cdot)$ is optimal in the chosen class
of controls, we must have
 \beq{diyi}\barray
W(s,x,p,\wdh {m}^\rho)\le V(s,x,p)+\frac{\rho}{3}.\earray\eeq

Note that for each given integer $\iota$, there is a measurable
function $F^\rho_\iota(\cdot)$ such that
$$\wdh {u}^\rho(t)=F^\rho_\iota(w_l(s),p(s),s\le \iota\delta,l\le
N)$$ on $[\iota\delta,\iota\delta+\delta)$. We next approximate
$F^\rho_\iota(\cdot)$ by a function that depends only on the sample
of $(w_l(\cdot),p(\cdot),l\le N)$ at a finite number of time points.
Let $\theta<\delta$ such that $\delta/ \theta$ is an integer.
Because the $\sigma-$ algebra determined by $\{{w_l(\nu
\theta),p(\nu\theta),\nu\theta\le \iota\delta, l\le N}\}$ increases
to the $\sigma$-algebra determined by $\{{w_l(s),p(s),s\le
\iota\delta, l\le N}\}$, the martingale convergence theorem implies
that for each $\delta, \iota$, there are measurable function
$F^{\rho,\theta}_\iota(\cdot)$, such that as $\theta \to 0$,
$$F^{\rho,\theta}_\iota(w_l(\nu\theta), p(\nu\theta), \nu\theta\le \iota\delta,l\le N)=u^{\rho,\theta}_\iota\to
\wdh {u}^\rho(\iota\delta)\ \hbox{ w.p.1.}$$ Here, we select
$F^{\rho, \theta}_\iota(\cdot)$ such that there are $N$ disjoint
hyper-rectangles that cover the range of its arguments and that
$F^{\rho, \theta}_\iota(\cdot)$ is constant on each hyper-rectangle.
Let $m^{\rho,\theta}(\cdot)$ denote the relaxed control
representation of the ordinary control $u^{\rho,\theta}(\cdot)$
which takes value $u^{\rho,\theta}_\iota$ on
$[\iota\delta,\iota\delta+\delta)$, and let
$(x^{\rho,\theta}(\cdot),p^{\rho,\theta}(\cdot))$ denote the
associated solution. Then for small enough $\theta$, we have
\beq{d2}\barray W(s,x,p,m^{\rho,\theta})\le W(s,x,p,\wdh
{m}^\rho)+\frac{\rho}{3}.\earray\eeq
 Next, we adapt
$F^{\rho,\theta}_\iota(\cdot)$ such that it can be applied to
$\{{\xi^{h_1,h_2}_n}\}$. Let $\bar{u}^{h_1,h_2}_n$ denote the
ordinary admissible control to be used for the approximation chain
$\{{\xi^{h_1,h_2}_n}\}$ defined in \eqref{fd}.

For $n$ such that $nh_2<\delta$, we can use any control. For
$\iota=1,2,\ldots$ and $n$ such that
$nh_2\in[\iota\delta,\iota\delta+\delta)$, use the control defined
by
$\bar{u}^{h_1,h_2}_n=F^{\rho,\theta}_\iota(w^{h_1,h_2}_l(\nu\theta),p^{h_2}(\nu\theta),\nu\theta
\le \iota\delta,l\le N)$. Recall that $\bar{m}^{h_1,h_2}(\cdot)$
denote the relaxed control representation of the continuous
interpolation of $\bar{u}^{h_1,h_2}_n$, then \bea
\ad(\xi^{h_1,h_2}(\cdot),\bar{m}^{h_1,h_2}(\cdot),
w^{h_1,h_2}_l(\cdot),
F^{\rho,\theta}_\iota(w^{h_1,h_2}_l(\nu\theta),p^{h_2}(\nu\theta),\nu\theta
\le \iota\delta,l\le N,\iota=0,1,2,\ldots  ))\\
\ad\to (x^{\rho,\theta}(\cdot),
m^{\rho,\theta}(\cdot),w_l(\cdot),F^{\rho,\theta}_\iota(w_l(\nu\theta),p(\nu\theta),
\nu\theta \le \iota\delta,l\le N,\iota=0,1,2,\ldots)).\eea Thus
$$W(s,x,p,\bar{m}^{h_1,h_2}) \le W(s,x,p,m^{\rho,
\theta})+\frac{\rho}{3}$$ Note that $$V^{h_1,h_2}(s,x,p)\le
W(s,x,p,\bar{m}^{h_1,h_2}).$$ Combing the above inequalities, we can
see $\limsup_{{h_1,h_2}\to 0}V^{h_1,h_2}(s,x,p)\le V(s,x,p)$ for
the chosen subsequence. By the tightness of $(\xi^{h_1,h_2}(\cdot),
p^{h_2}(\cdot),\bar{m}^{h_1,h_2}(\cdot))$ and arbitrary of $\rho$,
we get
$$\limsup_{{h_1,h_2}\to 0}V^{h_1,h_2}(s,x,p)\le V(s,x,p)$$ and thus
conclude the proof.\qed

\section{A Numerical Example}\label{sec:num}
\subsection{An Example}
In this section, we provide an example to demonstrate
our results.

\begin{exm}\label{7.1}
{\rm We consider a networked system with regime switching. There are
$2$ nodes in the system. One of the node has dynamic given by
$$dx_0(t)=r(t,\al(t))x_0(t)dt,$$ where $r(t,\al(t))=t+\al(t)$, the
other node follows the systems of SDEs
$$dx_1(t)=x_1(t)b(t,\al(t))dt+x_1(t)\sigma(t,\al(t))dw_1(t),$$ where
$b(t,\al(t))=1+t-\al(t)$, and $\sigma(t,\al(t))=\al(t)$. Observation
process is given by
$$dy(t)=g(\al(t))dt+dw_2(t),$$
with $g(1)=2$ and $g(2)=3$.
The Markov chain
$\al(\cdot)\in\{1,2\}$ is generated by the generator
$Q=\(\!\!\!\!\! \begin{array}{lrr} & -0.5 & 0,5 \\
& 0.5 & -0.5 \\ \end{array} \)$.

Our objective is to distribute proportions of the network flow to
each node so as to minimize the total variance at time $T$ subject to
$Ex(T)=\kappa$.
Our system $x(t)$ is $p^i(t)$ dependent and given by
$$dx(t)=[x(t)[(t+1)p^1(t)+(t+2)p^2(t)]-(p^1(t)+3p^2(t))u(t)]dt+u(t)[p^1(t)+2p^2(t)]dw_1(t).$$
To get the efficient frontier, note that on the one hand, $\kappa$
is given to us and we will choose a series of value for $\kappa$
starting from $[1,5.5]$. On the other hand, we need to know
$\lambda$, here we use simplex method to get the its value. Using
value iteration and policy iterations, we have the outline of the
procedure to find an improved values of $V$ as follows:

\beq{tran1}\barray \ad\hspace*{-1in} V^{h_1,h_2}(nh_2,x,p)\\
\aad\hspace*{-0.8in}  =
\min_{r\in \mathcal{U}^{h_1,h_2}}\sum_y  p^{h_1,h_2}((nh_2,x)(nh_2+h_2,y)|p,r)V^{h_1,h_2}(nh_2+h_2,y,p)\\
\aad\hspace*{-0.4in} +p^{h_1,h_2}(nh_2, y)|x,p,r)V^{h_1,h_2}(y,x,p)\\
\aad\hspace*{-0.4in} +\sum^m_{i=1}V^{h_1,h_2}(nh_2+h_2,x,p^i+h_1)\frac{\frac{1}{\sigma^2_0}[p^i(g(i)-\lbar\alpha)]^2h_2+2h_1(\sum_{j=1}^mq^{ji}p^j)^+h_2}{2h^2_1}\\
\aad \hspace*{-0.4in} +\sum^m_{i=1}V^{h_1,h_2}(nh_2+h_2,x,p^i-h_1)\frac{\frac{1}{\sigma^2_0}[p^i(g(i)-\lbar\alpha)]^2h_2+2h_1(\sum_{j=1}^mq^{ji}p^j)^-h_2}{2h^2_1}\\
\aad\hspace*{-0.4in} +\sum^m_{i=1}V^{h_1,h_2}(nh_2+h_2,x,p^i)[-\frac{\frac{1}{\sigma^2_0}[p^i(g(i)-\lbar\alpha)]^2h_2}{h^2_1}-\frac{h_2|\sum^m_{j=1}q^{ji}p^j|}{h_1}],\\
V^{h_1,h_2}(T,x,p)\ad=(x-\frac{1}{2})^2 \text{ for } x\not\in
[0,2].\earray\eeq

The corresponding control $u$ can be obtained as follows:
\beq{u}\barray \ad u^{h_1,h_2}(nh_2,x,p)\\
\aad\ =
\arg\min_{r\in \mathcal{U}^{h_1,h_2}}\sum_y p^{h_1,h_2}((nh_2,x)(nh_2+h_2,y)|p,r)V^{h_1,h_2}(nh_2+h_2,y,p)\\
\aad\quad +p^{h_1,h_2}(nh_2, y)|x,p,r)V^{h_1,h_2}(y,x,p)\\
\aad\quad +\sum^m_{i=1}V^{h_1,h_2}(nh_2+h_2,x,p^i+h_1)\frac{\frac{1}{\sigma^2_0}[p^i(g(i)-\lbar\alpha)]^2h_2+2h_1(\sum_{j=1}^mq^{ji}p^j)^+h_2}{2h^2_1}\\
\aad\quad +\sum^m_{i=1}V^{h_1,h_2}(nh_2+h_2,x,p^i-h_1)\frac{\frac{1}{\sigma^2_0}[p^i(g(i)-\lbar\alpha)]^2h_2+2h_1(\sum_{j=1}^mq^{ji}p^j)^-h_2}{2h^2_1}\\
\aad\quad +\sum^m_{i=1}V^{h_1,h_2}(nh_2+h_2,x,p^i)[-\frac{\frac{1}{\sigma^2_0}[p^i(g(i)-\lbar\alpha)]^2h_2}{h^2_1}-\frac{h_2|\sum^m_{j=1}q^{ji}p^j|}{h_1}].
\earray\eeq
\begin{figure}
\begin{center}
\includegraphics[width=80mm]{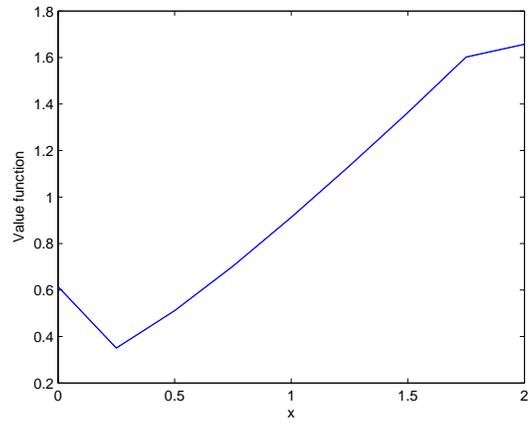}
\caption{Approximate value function with $h_1=0.25$  and $h_2=0.
001$ for fixed expectation}
\end{center}
\end{figure}

The value function is plotted in Figure 1, the corresponding control
in Figure 2, and the efficient frontier in Figure 3.

\begin{figure}
\begin{center}
\includegraphics[width=80mm]{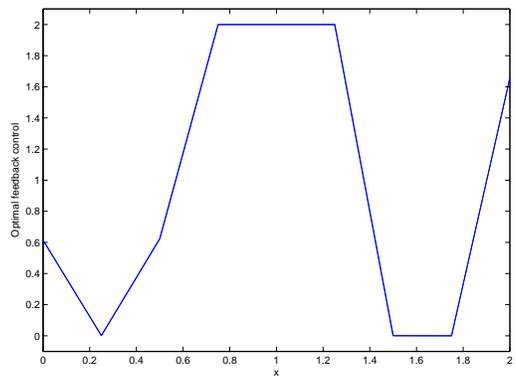}
\caption{Optimal feedback control with $h_1=0.25$
and $h_2=0.001$ for fixed expectation}
\end{center}
\end{figure}

\begin{figure}
\begin{center}
\includegraphics[width=80mm]{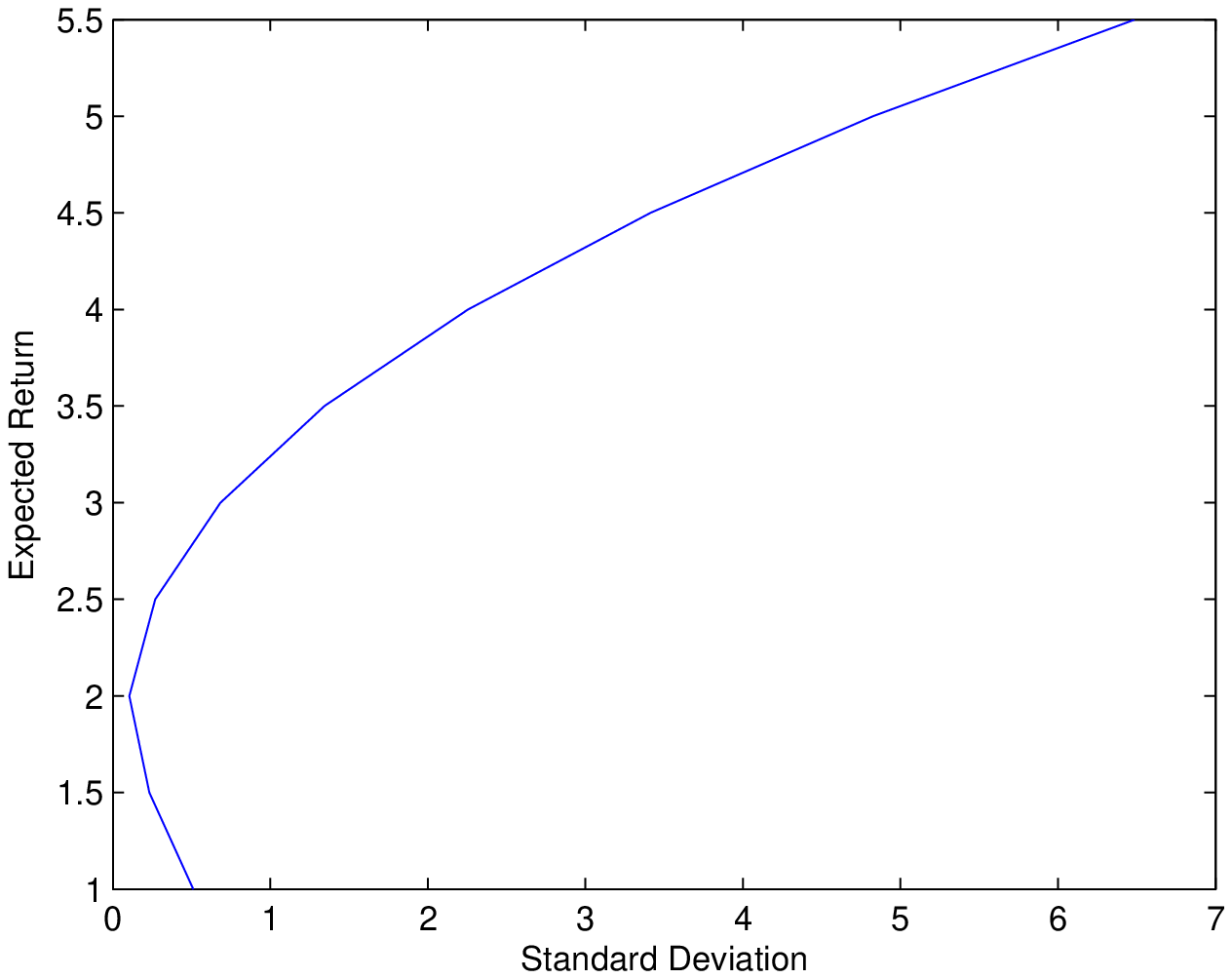}
\caption{Efficient frontier $h_1=0.25$  and $h_2=0.001$
when using simplex method to find out the optimal
$\lambda$}
\end{center}
\end{figure}

}
\end{exm}

\subsection{Further Remarks}
This paper developed a numerical approach for
a controlled switching diffusion system with a hidden Markov chain.
Using  Markov chain approximation techniques combined with the Wonham filtering,
a numerical scheme was developed.
In contrast to the existing work in the literature, we used
Markov chain approximation for the diffusion component and used a direct
discretization for the Wonham filter.
Our on-going effort will be directed to
use the approach developed in this work to treat certain networked systems
that involve platoon controls with wireless communications.


\begin{thebibliography}{99}

\setlength{\baselineskip}{0.12in}

{\small
\parskip=0pt





\bibitem{LZ}
S. Chen, X. Li, and X.Y. Zhou, stochastic linar quadratic
regulators with indefinite control weight costs, {\sl SIAM J. Control
Optim.}, {\bf 36} (1998), 1685--1702.


\bibitem{DB} F. Dufour and P. Bertrand, The filtering problem for
continuous time linear systems with Markovian switching
coefficients, {\sl Syst. Control Lett.}, {\bf 23} (1994),
453--461.


\bibitem{D}
D. Duffie and H. Richardson, Mean-variance hedging in continuous
time, {\sl Ann. Appl. Probab.}, {\bf 1} (1991), 1--15.

%
%
%

\bibitem{FN}
W.H. Fleming and M. Nisio,  On stochastic relaxed control for
partially observed diffusions, {\sl Nagoya Math. J.}, {\bf 93}
(1984), 71--108.









\bibitem{K}
H.J. Kushner, Numerical methods for stochastic control problems in
continuous time, {\sl SIAM J. Control Optim.}, {\bf 28} (1990),
999--1048.




\bibitem{K1}
H.J. Kushner, Consistency issues for numerical methods for variance
control with applications to optimization in finance, {\sl IEEE.
Trans. Automat. Control}, {\bf 44} (2000), 2283--2296.

\bibitem{KD} H.J. Kushner and P. Dupuis, {\it
Numerical Methods for Stochastic Control Problems in Continuous
Time}, Springer, New York, 2001.




\bibitem{LN}
D. Li and W.L. Ng, Optimal dynamic portfolio selection:
Multi-period mean-variance formulation, {\sl Math. Finance.}, {\bf
10} (2000), 387--406.





\bibitem{M}
H. Markowitz, Portfolio selection, {\sl J. Finance.}, {\bf 7}
(1952), 77--91.


\bibitem{S}
P.A. Samuelson, Lifetime portfolio selection by dynamic stochastic
programming, {\sl Rev. Econ. Statist.}, {\bf 51} (1969), 239--246.

\bibitem{P}
S.R. Pliska,  {\it Introduction to Mathematical Finance}, Basil
Blackwell, Malden, UK, 1997.


\bibitem{S1}
M. Schweizer, Approximation pricing and the variance optimal
martingale measure, {\sl Ann. Probab.}, {\bf 24} (1996), 206--236.

 \bibitem{SY} Q.S. Song, G. Yin, and Z. Zhang,
Numerical method for controlled regime-switching diffusions and
regime-switching jump diffusions,  {\sl Automatica}, {\bf 42} (2006),
1147--1157.


\bibitem{Wo}
W.M. Wonham,  Some applications of stochastic differential equations
to optimal nonlinear filtering, {\sl SIAM J. Control}, {\bf 2}
(1965), 347--369.


\bibitem{YYWZ} Z. Yang, G. Yin,  L,Y, Wang, and H. Zhang,
 Near-Optimal mean-variance controls
under two-time-scale formulations and applications,
 {\sl Stochastics}, {\bf 85} (2013), 723--741.

\bibitem{DZ}
D. Yao, Q. Zhang, X.Y. Zhou,
A regime switching model for European options,
in {\it Stochastic Process, Optimization, and Control Theory},
 H. Yan, G. Yin and Q. Zhang Eds., 281--300, Springer, 2006.


\bibitem{YinZ} G. Yin and Q. Zhang,
{\it Discrete-Time Markov chains: Two-time-scale methods and
application}, Springer, New York, 2005.

\bibitem{Y2}
G. Yin and X.Y. Zhou, Markowitz mean-variance portfolio selection
with regime switching: from discrete-time models to their
continuous-time limits, {\sl IEEE Trans. Automatic
Control.}, {\bf 49} (2004), 349--360.


\bibitem{YZ}
G. Yin and C. Zhu, {\it Hybrid Switching Diffusions.} New York:
Springer, 2010.


\bibitem{JZ}
J. Yong and X.Y. Zhou, {\it Stochastic controls: Hamiltonian Sytems
and HJB Equations.} New York: Springer, 1999.



\bibitem{LZY}
L. Yu, Q. Zhang and G. Yin, Asset allocation for regime switching
market models under partial observation, to appear in {\sl Dynamic Systems and
Applications.}

\bibitem{Z}
Q. Zhang, Stock trading: an optional selling rule, {\sl SIAM
J.Control Optim.}, {\bf 40} (2001), 64--87.


\bibitem{ZL} X.Y. Zhou and D. Li,
Continuous time mean variance portfolio selection: A stochastic LQ
Framework, {\sl Appl Math Optim}, {\bf 42} (2000), 19--33.

\bibitem{ZY} X.Y. Zhou and G. Yin,
Markowitz mean-variance portfolio selection with regime switching: A
continuous time model, {\sl SIAM J.Control Optim.}, {\bf 42} (2003),
1466--1482.

 }
\end{thebibliography}
\end{document}